\documentclass[conference]{IEEEtran}
\IEEEoverridecommandlockouts
\usepackage{cite}
\usepackage{amsmath,amssymb,amsfonts}
\usepackage{algorithmic}
\usepackage{graphicx}
\usepackage{textcomp}
\usepackage{xcolor}
\usepackage{amsbsy}
\usepackage{amsmath,amsfonts}
\usepackage{algorithmic}
\usepackage{graphicx}
\usepackage{textcomp}
\usepackage{xspace}
\usepackage{comment}
\usepackage{enumitem}
\usepackage{mdframed}
\usepackage{url}
\usepackage{todonotes}
\usepackage{color, colortbl}
\usepackage{listings}
\usepackage{multirow}
\usepackage[hidelinks]{hyperref}
\usepackage{booktabs}
\usepackage{adjustbox}

\def\BibTeX{{\rm B\kern-.05em{\sc i\kern-.025em b}\kern-.08em
    T\kern-.1667em\lower.7ex\hbox{E}\kern-.125emX}}
\begin{document}

\title{Testing GPU Numerics: Finding Numerical Differences Between NVIDIA and AMD GPUs
\thanks{This 
work was performed under
the	auspices of the	U.S. Department of Energy by Lawrence Livermore	
National Laboratory under Contract DE-AC52-07NA27344 (LLNL-CONF-868447).}
}


\author{
    \IEEEauthorblockN{Anwar Hossain Zahid}
    \IEEEauthorblockA{
        \textit{Department of Computer Science} \\
        Iowa State University \\
        Ames, IA \\
        ahzahid@iastate.edu
    }
    \and
    \IEEEauthorblockN{Ignacio Laguna}
    \IEEEauthorblockA{
        \textit{Center for Applied Scientific Computing} \\
        Lawrence Livermore National Laboratory \\
        Livermore, CA \\
        ilaguna@llnl.gov
    }
    \and
    \IEEEauthorblockN{Wei Le}
    \IEEEauthorblockA{
        \textit{Department of Computer Science} \\
        Iowa State University \\
        Ames, IA \\
        weile@iastate.edu
    }
}


\maketitle


\definecolor{dkgreen}{rgb}{0,0.6,0}
\definecolor{darkred}{rgb}{0.3,0.1,0.1}
\definecolor{gray}{rgb}{0.5,0.5,0.5}
\definecolor{mauve}{rgb}{0.58,0,0.82}
\definecolor{light-gray}{gray}{0.9}
\definecolor{blue}{rgb}{0,0,0.75}
\definecolor{peach}{rgb}{1.0,0.91,0.84} 
\newcommand{\wei}[1]{\textcolor{blue}{[Wei]: #1}}
\lstset{ %
  language=C,
  basicstyle=\scriptsize\ttfamily,
  numbers=left,
  numberstyle=\scriptsize\color{gray},  
  numbersep=5pt,                  
  backgroundcolor=\color{white},
  showspaces=false,               
  showstringspaces=false,         
  showtabs=false,
  frame=single,                   
  rulecolor=\color{black},
  tabsize=2,                      
  captionpos=b,                   
  breaklines=true,                
  breakatwhitespace=false,
  keywordstyle=\color{blue},          
  commentstyle=\color{dkgreen},       
  stringstyle=\color{mauve},         
  escapeinside={\%*}{*)},
  xleftmargin=4.0ex,
  morekeywords={for,each,between,can,reach,in,is,Sort,Print,From}
}

\providecommand{\en}[1]{\ensuremath{\text{E{#1}}  }}

\newcommand{\subheader}[1]{\noindent \textbf{#1}}

\newcommand{\squeezeup}{}

\hyphenation{allo-ca-tion scale-de-pen-dent}

\renewcommand{\ttdefault}{cmtt}

\newlist{rqs}{enumerate}{1}
\setlist[rqs,1]{label={\bfseries Q\arabic*},align=left,wide, labelwidth=!, 
labelindent=0pt}

\newcommand{\note}[1]{{\color{red} NOTE: #1}}

\newcommand{\infinity}{\texttt{INF}\xspace}
\newcommand{\positiveinf}{\texttt{INF+}\xspace}
\newcommand{\negativeinf}{\texttt{INF-}\xspace}
\newcommand{\nan}{\texttt{NaN}\xspace}
\newcommand{\nans}{\texttt{NaNs}\xspace}
\newcommand{\hipify}{\texttt{HIPIFY}\xspace}
\newcommand{\HIPIFY}{\texttt{HIPIFY}\xspace}
\newcommand{\varity}{\texttt{Varity}\xspace}

\pagestyle{plain}

\begin{abstract}
As scientific codes are ported between GPU platforms,
continuous testing is required to ensure numerical robustness and 
identify potential numerical differences between platforms. 
Compiler-induced numerical differences can occur when 
a program is compiled and run on different GPUs and compilers,
and the numerical outcomes are different for the same input.
We present a study of compiler-induced numerical differences
between NVIDIA and AMD GPUs, two widely used GPUs in HPC clusters.
Our approach uses a random program generator (\varity) to generate
thousands of short numerical tests in CUDA and HIP, and 
their inputs; then, we use differential testing to check if the  program produced a numerical 
inconsistency when run on NVIDIA and AMD GPUs, using the same
compiler optimization level.
We also use the AMD's \hipify tool to convert 
CUDA tests into HIP tests and test if there are numerical inconsistencies induced by \hipify.
%
%
In our study, we generated more than 600,000 tests and 
found subtle numerical differences occurring between 
the two classes of GPUs. We found that some of the differences come
from (1) math library calls, (2) differences in floating-point 
precision (FP64 versus FP32),
and (3) converting code to HIP with \hipify.
\end{abstract}

\begin{IEEEkeywords}
Differential testing, HIP, CUDA, random program generation,
floating-point.
\end{IEEEkeywords}

\section{Introduction}

Testing scientific software on different platforms is crucial 
for developers to ensure the consistency, accuracy, and reliability of
the numerical results. As scientific 
software is ported to heterogeneous platforms involving 
different classes of GPUs, it is essential to 
test the programming environment in such platforms to understand
when numerical differences or reproducibility issues emerge.
Previous studies have shown that scientific applications
can output numerical results that differ from each other when
run in different heterogeneous systems~\cite{ahn2021keeping,miao2023expression,guo2020pliner}, depending
on several factors, such as the compilers used to emit code,
the compiler optimizations used, and how the GPU
implements floating-point arithmetic.

While NVIDIA GPUs are the most used in HPC clusters,
AMD GPUs have been the choice for building recent supercomputers.
Exascale computing platforms from the US Department of Energy (DOE), 
such as Frontier and El Capitan use AMD GPUs (MI250 and MI300 GPUs, 
respectively). AMD GPUs are also widely used in several areas of 
scientific research, including machine learning, climate research, and genomics. As scientists port codes from NVIDIA GPUs to AMD GPUs,
or vice versa, it is crucial to understand the numerical differences
and numerical reproducibility challenges that developers can observe
in these platforms.

Previous work has studied numerical portability between various GPUs,
including AMD GPUs. For example, authors in~\cite{innocente2021accuracy} 
target math function results relative to earlier versions of GPUs.
Li et al.~\cite{li2024fttn} create tests to compare numerical differences
and capabilities in NVIDIA and AMD GPUs with a focus on matrix accelerators
and tensor cores. While these studies reveal interesting differences
between these classes of GPUs, they have drawbacks: (1) tests
are manually generated, which limits the number of floating-point arithmetic
expressions that are tested; (2) they do not consider compiler optimizations
between the GPU compilers---studies have shown that compiler optimizations
can be the source of significant differences between GPU platforms~\cite{guo2020pliner,ahn2021keeping,bentley2019multi}.

\begin{figure*}[th!]
    \centering
    \includegraphics[width=0.85\linewidth]{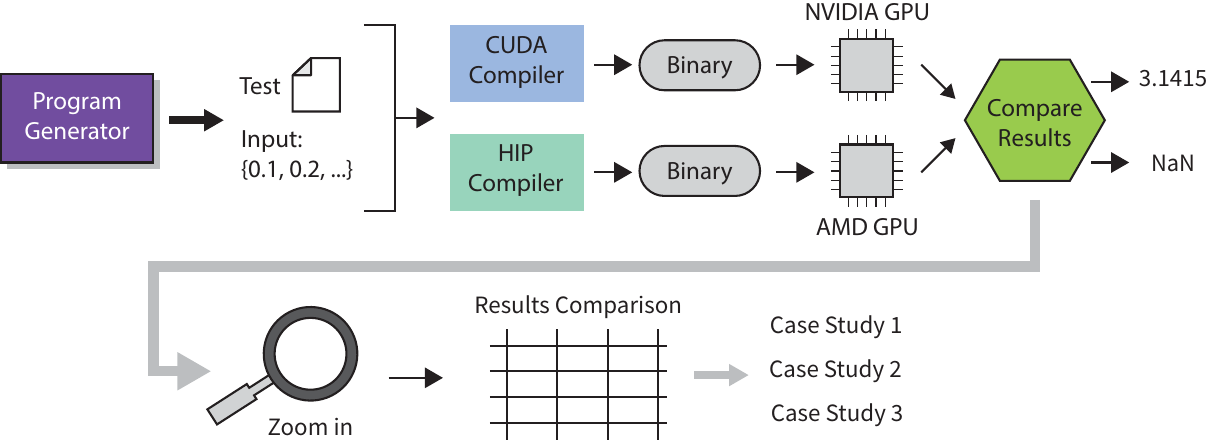}
    \caption{Overview of testing approach via random program generation for both GPUs (NVIDIA and AMD).}
    \label{fig:overview}
\end{figure*}

Random program generation has been used before in the \varity~\cite{varity}
framework to test numerical differences that are induced by compilers
between NVIDIA GPUs and CPUs. The approach generates thousands of random
floating-point programs and their numerical inputs, compiles them with multiple
optimization flags, runs them on NVIDIA GPUs and CPUs, and identifies cases
that induce different numerical results. The work in~\cite{varity}, however, 
did not include tests for HIP, the programming interface for AMD GPUs.

In this paper, we use the approach of random program generation and differential testing
to expose \textit{compiler-induced numerical inconsistencies} 
between NVIDIA and AMD GPUs. These inconsistencies arise when a given test program
$P$ and its input $I$, compiled with the same optimization level and run
on two different platforms, produce different numerical results. For example,
consider a test program $P$ with inputs $I=\{0.1, 0.2, 0.3\}$. Suppose
we compile $P$ with the HIP compiler \texttt{hipcc}, which produces the binary
$P_{\textit{HIP}}$. Compiling the CUDA version of $P$ with \texttt{nvcc} (the NVIDIA compiler)
produces the binary $P_{\textit{NVCC}}$. Examples of compiler-induced inconsistencies occur
when both $P_{\textit{HIP}}$ and $P_{\textit{NVCC}}$ 
are run with the same input $I$, but they produce
different outputs, e.g., $3.1415$ versus $3.999$, or $3.1415$ versus a 
not-a-number (\texttt{NaN}).

\subheader{Summary of Contributions.}
We present a test-guided study about the compiler-induced
numerical differences between NVIDIA and AMD GPUs.
We use a random program generation approach to generate \textbf{thousands} of 
floating-point arithmetic tests that expose the numerical differences developers
can experience when running numerical codes on these GPUs.
When such numerical differences are observed, the framework provides the user with
a \textit{small test} with the inputs that induce such inconsistency.
Such small tests have several benefits for testing HPC systems: (a)
the tests are easier to analyze than a large scientific
application---they are composed of a single kernel and a single set
of numerical inputs; (b) the tests can be provided to vendors 
for further investigation---they are self-contained; (c) they can be
reused to test newer systems, potentially during acceptance testing.

Our study uses the \varity~\cite{varity} framework to generate
such tests. While \varity originally generated CUDA programs,
it did not generate HIP tests. We have extended it to generate 
HIP tests and tested it with the recent ROCm compilers. 
We perform a large-scale testing campaign and run more than {600,000}
tests in these two classes of GPUs. We present several case studies showing how these
numerical differences can emerge when using different compiler optimizations.

In addition to generating HIP tests directly 
through our \textit{extended Varity framework} (hereafter referred to as \varity), we use AMD's \texttt{HIPIFY}\footnote{\url{https://rocm.docs.amd.com/projects/HIPIFY/en/latest/}}
tool to convert CUDA source code into HIP. \texttt{HIPIFY} facilitates the 
porting process between NVIDIA and AMD platforms by automatically 
translating CUDA applications into HIP. We ran 123,750 tests on \texttt{HIPIFY}-converted
programs to compare the numerical results with those produced by 
the HIP tests generated by \varity, aiming to 
identify any discrepancies introduced by the \texttt{HIPIFY} conversion.

In summary, our contributions are the following:
\begin{itemize}
    \item We present an approach to test the numerical differences induced
    by compiler optimizations between two widely used GPUs in HPC clusters,
    NVIDIA and AMD GPUs. Our approach extends the \varity framework
    and adds support for testing HIP programs.
    \item We evaluate the approach in two production clusters at the Lawrence Livermore
    National Laboratory, one with NVIDIA V100 GPUs (Lassen)
    and the other one with AMD MI250 GPUs (Tioga). We present a detailed evaluation 
    and results for different levels of optimization on these platforms.
    \item We use the AMD's \texttt{HIPIFY} tool to convert CUDA source code 
    into HIP, identifying and analyzing discrepancies introduced by the automatic code translation.
    \item We present several case studies that reveal numerical 
    inconsistencies on these GPUs, highlighting discrepancies caused by 
    math functions (e.g., \texttt{fmod} and \texttt{ceil}), the effects 
    of different optimization levels, and variations in handling 
    special values such as \texttt{NaN} and \texttt{Infinity}.
\end{itemize}

\section{Background and Overview}

In this section, we present background information needed
to understand our approach and an overview of the approach.

\subsection{Compiler-Induced Numerical Differences}

Compilers used in heterogeneous HPC platforms (e.g., \texttt{clang},
\texttt{nvcc}, \texttt{hipcc}) provide several levels of optimization flags from
\texttt{-O0} to \texttt{-O3}. With higher optimization levels, program
performance can be improved at the cost of
potentially generating non-compliant IEEE 754 floating-point
code. There are optimization flags that explicitly violate the IEEE
754 standard, but can offer significant speedups. 
We say that a \textit{compiler-induced numerical difference} (or inconsistency)
occurs when a program that is tested on two different compilers (e.g., a compiler for GPU 1
and another compiler for GPU 2) produced different numerical results,
even when the same optimization level is used, e.g., \texttt{-O0}.

\begin{table}[th]
 \centering
\renewcommand{\arraystretch}{1.0}
\caption{Inconsistencies in BT.S.}
\scalebox{1.2}{
\begin{tabular}{llc}
    \toprule
Compiler Options     & Runtime & Error \\ 
\midrule
nvcc -O0                  & 0.104s  & 6.98176E-13    \\
nvcc -O3 -use\_fast\_math & 0.052s  & 9.73738E-13    \\
clang -O0                 & 0.349s  & 8.32928E-13    \\
clang -O3 -ffast-math     & 0.059s  & 3.50905E-12    \\ 
\bottomrule
\end{tabular}
}
\label{table:inc_example}
\end{table}

Table~\ref{table:inc_example} (taken from~\cite{miao2023expression})
shows an example of a numerical
inconsistency between an NVIDIA GPU and a CPU in the BT NAS program---the compiler used
in the NVIDIA GPU is \texttt{nvcc} and for the CPU \texttt{clang}.
The Table shows the program runtime and the maximum relative error for each compiler
and optimization flag combination. Using \texttt{-O3 -use\_fast\_math}
with \texttt{nvcc} yields 100\% speedup compared to \texttt{-O0}, but at the
cost of the error being 39\% larger. The performance and error with Clang
are generally worse, with the largest error in \texttt{clang -O3
-ffast-math} being 402\% larger than \texttt{nvcc -O0}.

In real-world scientific applications, such compiler-induced numerical
inconsistencies are not uncommon. They can happen when migrating
software to other hardware/software platforms, switching applications
to a new compiler, or just using more aggressive optimization flags
for compilation. These inconsistencies may cause major software
failures that take a tremendous amount of effort to identify and
resolve (see~\cite{ahn2021keeping}).

\subsection{IEEE 754 Exceptions and Optimizations}

A floating-point number has the form
\begin{equation}
x = \pm m \times \beta^e \label{lab:fpformat}
\end{equation}
where the sign, the mantissa $m$, the exponent $e$ can be stored in memory or a 
register. We assume $\beta=2$, since this is the most used format for representing 
floating-point numbers.
The IEEE 754 Standard defines five classes of \textit{exceptions}, 
that can result from arithmetic operations. 
Table~\ref{table:fp_events} shows these five 
events. 

\begin{table}
	\caption{IEEE 754 Standard exceptions.}
	\begin{center}
		\scalebox{0.95}{
			\begin{tabular}{lll}
				& \textbf{Event} & \textbf{Description} \\ \hline
				\multirow{5}{2em}{IEEE 754} & \texttt{Inexact} & Result is produced after 
				rounding \\
				& \texttt{Underflow} & Result could not be represented as normal \\
				& \texttt{Overflow} & Result did not fit and it is an infinity \\
				& \texttt{DivideByZero} & Divide-by-zero operation \\
				& \texttt{Invalid} & Operation operand is not a number (NaN) \\
				\hline
			\end{tabular}
		}
	\end{center}
	\label{table:fp_events}
	\vspace{-0.3cm}
\end{table}

The IEEE 754 Standard defines five classes of \textit{exceptions}, 
that can result from arithmetic operations. Table~\ref{table:fp_events} shows these five 
events. When one of these events occurs, the floating-point unit can set a status register 
specifying which event occurred. Existing frameworks for CPU analysis, such 
as~\cite{dinda2020spying} read these registers to detect the occurrence of such events. 
Additionally, with the help of the compiler and system routines, they raise a 
floating-point exception signal (e.g., \texttt{SIGFPE}) when these events 
occur in the CPU. Unlike CPUs, NVIDIA GPUs have no mechanism to detect floating-point 
exceptions, set a status register, or raise a signal when an exception occurs.

\subsubsection{\textbf{Exceptional Quantities}}
Except for the \texttt{Inexact} exception in Table~\ref{table:fp_events}, the rest of the 
events will result in either a \nan (not a number), \infinity (infinity, positive or 
negative), or a subnormal number (i.e., a number smaller than a normal floating-point 
number but that is 
not zero). More specifically, \texttt{Overflow} and \texttt{DivideByZero} result in \infinity, 
and \texttt{Invalid} result in \nan. \texttt{Underflow} can result in zero or a subnormal 
number---in our case, we are interested in underflows that result in subnormal numbers.
The \texttt{Inexact} event results in a rounding operation; however, this occurs frequently 
in numerical programs, and it is usually of no interest to programmers. In summary, our 
goal is to identify inputs that produce any of these cases: \nan, \positiveinf, 
\negativeinf, or subnormal 
number quantities (positive or negative).

\subsubsection{\textbf{Subnormal Numbers}}
Note that while subnormal numbers can represent specific real number quantities, they are 
often dangerous for several reasons. First, they indicate that computational
results are becoming too small to be represented in the current precision. Second, if they 
propagate to 
the denominator of a division, the result can produce \infinity. For example, $\frac{1}{x}$, 
where $x=$1e-309 (an FP64 subnormal number), produces \infinity\footnote{This behavior 
may depend on the platform, compiler, and optimizations applied to the code.}. Third, 
when 
subnormal 
numbers are combined with compiler optimizations, and they can cause reproducibility 
issues~\cite{guo2020pliner}. Therefore, it is crucial to mitigate them and understand 
when they occur, as well as \nan and \infinity.

\subsection{Approach Overview}
Figure~\ref{fig:overview} shows an overview of the approach.
The program generator first generates the source code of tests and inputs
into the tests. A file test is generated with the extension \texttt{.cu} for CUDA and
\texttt{.hip} for HIP. Then the tests are compiled with the corresponding
compilers (\texttt{nvcc} and \texttt{hipcc}) using the same optimization label, 
which produces two binaries.
The binaries are executed on the NVIDIA and AMD GPUs using the same input.
We then compare the numerical results of the tests and find anomalies,
i.e., differences.
After all tests are executed, meta-data is saved in a JSON file with all
the results. We then analyze them and identify several case studies 
(see Section~\ref{sec:cases}).

\section{Approach}

In this section, we present the details of our
approach.
First, we give a high-level overview of \varity~\cite{varity},
the approach we use to generate tests. Next, we explain
our implementation of HIP test generation in the framework.

\subsection{Varity's Framework}
\label{sec:varity}
Here, we present a high-level description of the
\varity framework. More details can be found
in the paper~\cite{varity}. 
\varity generates random programs that expose a 
wide range of floating-point arithmetic operations and 
other structures encountered in scientific codes,
such as \texttt{for} loops and \texttt{if} conditions.
\varity also generates random floating-point
inputs for the programs. \varity was originally designed 
to generate C/C++ tests and CUDA tests.
In this work, we extend it to support the generation of HIP tests.

\subsubsection{\textbf{Grammar to Define Possible Tests}}
A grammar is defined to specify the set of possible tests (or programs)
that the tool can generate. Due to space limitations, we do not show
the grammar in his paper, but it can be found in~\cite{varity}.
\varity's grammar considers the most important 
aspects of HPC programs and uses the characteristics of programs 
that could (most likely) affect how floating-point code is 
generated and executed. The characteristics of the programs
that \varity generates are shown in Table~\ref{tab:programs}.

\begin{table}[]
\centering
\caption{Characteristics of the Random Programs}

\scalebox{1.15}{
\begin{tabular}{p{1.5cm}|p{5cm}} \hline

\textbf{Floating-Point Types} & Variables 
using single and double floating-point precision 
(i.e., \texttt{float} and \texttt{double}). \\ \hline

\textbf{Arithmetic Expressions} & Arithmetic expressions can use any 
operator in \{\texttt{+}, \texttt{-}, \texttt{*}, \texttt{/}\}, 
parenthesis ``( )'', and functions from the \texttt{C math} library. 
The grammar also allows boolean expressions. \\ \hline

\textbf{Loops} & \texttt{for} 
loops with multiple levels of nesting. 
We can generate loop sets $L_1 > L_2 > L_3 > \ldots > L_N$, where 
$L_1$ encloses $L_2$, $L_2$ encloses $L_3$, 
and so on up to $L_N$, where $N$ is defined by the user. \\ \hline

\textbf{Conditions} & \texttt{if} conditions, 
which can be true or false based on a boolean expression. \\ \hline

\textbf{Variables} & Programs can contain temporal 
floating-point variables. Variables can store arrays 
or single values. \\ \hline

\end{tabular}
}
\label{tab:programs}
\end{table}

\subsection{Program Output}
All operations are enclosed in a kernel function named \texttt{compute}. 
The kernel function does not return anything; instead, it computes a 
floating-point value and stores it in the 
\texttt{comp} variable. The value of \texttt{comp} is printed in standard output.

Note that, in addition to the \texttt{comp} kernel function, 
the generators generate a \texttt{main()} function and code to 
allocate and initialize arrays (if arrays are used in the test program). 
For simplicity, we do not present this in the grammar. The 
\texttt{main()} function reads the program inputs and copies them to the 
\texttt{comp} kernel function parameters before calling the kernel function.
Figure~\ref{fig:example_test} shows a sample random test.

\begin{figure} 
    \begin{lstlisting}[
        xleftmargin=0.15in,
        linewidth=0.48\textwidth, 
        basicstyle=\scriptsize\ttfamily
        ] 
__global__
void compute(double comp, int var_1, double var_2,
  double var_3, double var_4, double var_5, 
  double var_6, double var_7, double var_8) {
  if (comp == -1.3857E-36 + var_2) {
    double tmp_1 = +1.3305E12 / var_3;
    comp += -1.7744E-2 * tmp_1;
    comp += cos(var_4 - +1.4014E2 * (var_5 + var_6 * var_7));
    for (int i=0; i < var_1; ++i) {
      comp -= sqrt(var_8 + -1.7976E3);
    }
  }
  printf("%.17g\n", comp);
}
\end{lstlisting}
\caption{Example of a simple test random program in FP64 precision.}
\label{fig:example_test}
\end{figure}

\subsection{FP32 and FP64 Support}
Our approach supports the generation of test programs in FP32 and FP64
precision. This is crucial because some compiler optimizations
may behave differently depending on the precision of arithmetic
operations. There is a configuration option that controls the types
of variables and functions used in the generation phase. When using
FP64, all types are \texttt{double}. When using FP32, all types
are \texttt{float}; this case includes using the FP32 math functions
where function names and constants end with \texttt{f}, e.g., \texttt{cosf()}, 
instead of \texttt{cos()}, and $1.23F$ instead of $1.23$.

\subsection{HIP Extensions}
To support HIP, we take advantage of the fact that HIP is considered
to be a subset of CUDA; thus, several of the CUDA constructs can be
reused in HIP---for example, a GPU kernel is declared as 
\texttt{\_\_global\_\_} in both APIs. There are, however, a few differences
in the API that we considered in our \varity implementation,
for example, the kernel launch API in CUDA uses \texttt{<<< >>>},
while HIP uses \texttt{hipLaunchKernel}.

\subheader{Compiler Matching.}
Generating random tests for a given platform involves deciding which
compiler to use for CUDA or HIP. Compiler matching is done automatically
depending on the program extensions---a random test file ends with \texttt{.cu}
is automatically compiled with \texttt{nvcc}, while HIP files are compiled
with \texttt{hipcc}, compiler driver utility from AMD ROCm.

\subheader{Fast Math Optimizations.}
The \texttt{-ffast-math} flag is a compiler option that 
enables a suite of optimizations designed to improve 
the performance of floating-point arithmetic 
by making several assumptions about the mathematical 
operations, such as no NaNs or Infinities, relaxed precision, and others. This flag includes optimizations like \texttt{-fno-math-errno}, 
\texttt{-funsafe-math-optimizations}, \texttt{-fno-trapping-math}, 
\texttt{-fassociative-math}, \texttt{-freciprocal-math}, \texttt{-fno-signed-zeros}, \texttt{-fno-rounding-math}, and \texttt{-ffinite-math-only}. 
While these optimizations work effectively in CUDA, 
they create issues in HIP when dealing with NaNs and Infinities, 
particularly due to the \texttt{-ffinite-math-only} flag. 
As \varity generates tests that may 
produce NaNs or infinities, the \texttt{-ffast-math} flag in HIPCC 
leads to errors. Therefore, the ROCm developers recommend using the \texttt{-DHIP\_FAST\_MATH}\footnote{\url{https://github.com/ROCm/HIP/issues/28}} flag instead, which provides similar performance benefits without causing the same issues. 
In summary, while CUDA handles the \texttt{-ffast-math} flag without problems, we are bound to use \texttt{-DHIP\_FAST\_MATH} in HIP to avoid issues with special floating-point values.

\subsection{Between-Platform Comparisons}
\label{subsec:between-platform-comparisons}
GPUs from different vendors are usually installed on different HPC 
clusters, a cluster with NVIDIA GPUs typically does not have AMD GPUs
and vice versa. Comparing numerical results between these GPUs require
running the same random tests in both clusters. Our testing approach
supports this by storing metadata of the tests generated in cluster $C_1$,
which is then used in cluster $C_2$ to run exactly \textit{the same} tests
and inputs.

\begin{figure}
    \centering
    \includegraphics[width=0.75\linewidth]{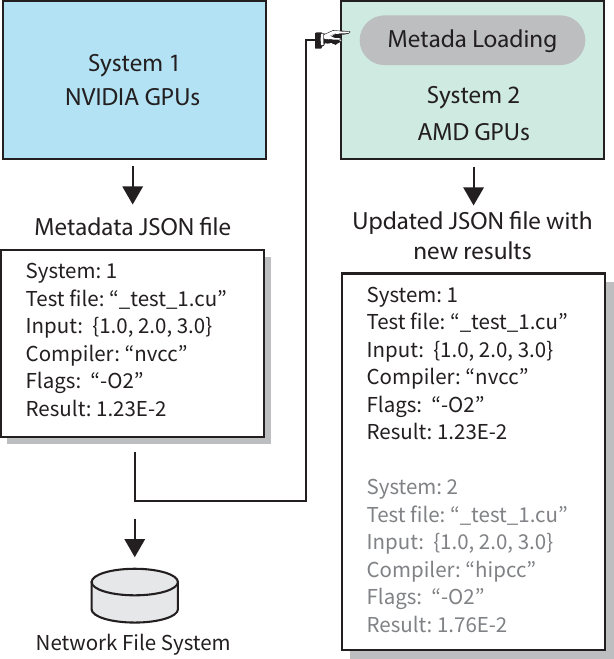}
    \caption{Process to perform between-platform comparisons.}
    \label{fig:metadata}
\end{figure}

Figure~\ref{fig:metadata} shows the process of between-platform testing.
We first run in system $C_1$. After all experiments are run,
\varity saves a JSON metadata file containing details about the
tests, inputs, compilers used, and results. Then, this metadata, along with the
tests, is transferred to system $C_2$. In system $C_2$, we load the metadata,
locate the tests we used in $C_1$, recompile them with the appropriate compiler,
and run the experiments again using different GPUs. We then save a new
metadata file that contains all the results. We use this new JSON file
to analyze the findings and search for cases with numerical inconsistencies.

\subsection{Enabling HIPIFY}
\hipify is a set of tools that developers can use to automatically 
translate CUDA source code into HIP source code. 
Since application developers may use \hipify to translate CUDA applications
into HIP applications, it is important to understand if numerical differences
occur in this approach.
We implement two approaches to generate HIP tests. The first approach
simply creates random tests using the random test generator we use for
generating CUDA tests. The second approach first generates a CUDA test
and then uses \hipify to translate the randomly generated CUDA 
source code into HIP source code.
\section{Evaluation}
\label{sec:evaluation}
In this section, we evaluate our approach on two high-performance computing systems, each employing NVIDIA and AMD GPUs. We begin by summarizing the numerical inconsistencies identified in our experiments, followed by an in-depth analysis of three case studies that provide insights into the sources of these discrepancies.

Our evaluation addresses the following research questions:
\begin{rqs}
\item How effective is our approach in detecting compiler-induced numerical inconsistencies between NVIDIA and AMD GPUs?
\item What specific classes of inconsistencies are revealed through our testing?
\item Can we determine the root causes of these inconsistencies?
\end{rqs}

\subsection{Systems and Software}

\subsubsection{System 1: Lassen with NVIDIA GPUs}
We conducted our experiments on the Lassen system at Lawrence Livermore National Laboratory (LLNL). Lassen is equipped with NVIDIA V100 GPUs and is a smaller-scale version of the classified Sierra system, offering a peak performance of 23 petaflops compared to Sierra's 125 petaflops. Notably, Lassen was ranked \#10 on the June 2019 Top500 list. The system comprises IBM Power9 CPUs with 44 cores per node, totaling 34,848 cores, alongside 3,168 NVIDIA V100 GPUs. We utilized CUDA version 12.2.2, the latest available on Lassen, and ran experiments using Python 3.9.12. The built-in Python libraries were sufficient for our needs, eliminating the requirement for a separate Python environment. The operating system used was Red Hat Enterprise Linux (RHEL).

\subsubsection{System 2: Tioga with AMD GPUs}
For experiments involving AMD hardware, we used the Tioga system at LLNL, which features 128 AMD MI-250X GPUs. Tioga is powered by AMD Trento CPUs, with 64 cores per node, amounting to a total of 2,048 cores. We used ROCm version 6.1.2, AMD's open-source platform for GPU programming, to run our experiments. The system operates on TOSS 4, a customized version of the Linux operating system. As with Lassen, all the tests on Tioga were conducted using Python 3.9.12, and the default Python libraries were adequate for our experiment.

\subsection{Experiments Configuration}

We executed a total of \textbf{652,600} experimental instances, 
each consisting of a program and input combination. 
We tested these instances under two settings: double precision (FP64) 
and single precision (FP32). 
We also used the \HIPIFY tool to convert FP64 tests from CUDA to HIP.

Due to resource constraints, we divided the tests 
into multiple batches, executed each batch separately, 
and then compiled the results into a comprehensive 
dataset, as described in Subsection~\ref{subsec:between-platform-comparisons}. We evaluated each generated 
code across \textbf{five optimization levels}: 
\texttt{O0}, \texttt{O0}, \texttt{O2}, \texttt{O3}, 
and \texttt{O3} with the \texttt{-ffast-math} flag. 
\texttt{O0} represents no optimization, while \texttt{O1} through \texttt{O3} denote progressively more aggressive optimization levels.
The \texttt{-ffast-math} flag further increases optimization aggressiveness by prioritizing speed over numerical precision, potentially sacrificing accuracy for performance gains.

We identified four possible outcomes from any test: 
\texttt{NaN} (Not a Number), \texttt{Inf} (Infinity), 
\texttt{Zero}, and \texttt{Number}\footnote{The term ``Number" refers 
to a non-zero real-valued floating-point number.}.
We categorized the observed discrepancies into \textbf{seven distinct types} 
based on the results from the \texttt{nvcc} and \texttt{hipcc} compilers: 
\texttt{NaN vs. Inf}, \texttt{NaN vs. Zero}, \texttt{NaN vs. Number}, \texttt{Inf vs. Zero}, \texttt{Inf vs. Number}, \texttt{Zero vs. Number}, and \texttt{Number vs. Number}. 

We make this distinction between \texttt{Zero} and \texttt{non-zero 
floating point numbers} to emphasize how computations involving 
near-subnormal or near-infinite values can lead to zero 
results due to rounding or precision errors. 
This distinction is important because extreme values often 
result in zero. 
However, we excluded inconsistencies such as \texttt{-NaN vs. +NaN}, 
\texttt{-Inf vs. +Inf}, and \texttt{-Zero vs. +Zero}, 
as these do not represent true numerical differences. 


\subsection{Experimental Results}
Table~\ref{tab:sum_exp_res} summarizes the experimental results, 
detailing key metrics for FP64, \HIPIFY-converted FP64, 
and FP32 tests across various optimization levels. We conducted a total 
of \texttt{652,600} runs: \texttt{247,500} each for FP64 
and \HIPIFY-converted FP64 tests, and \texttt{157,600} for FP32 tests. 
We observe discrepancies in \textbf{0.98\%} of the FP64 runs, 
\textbf{1.10\%} of \HIPIFY-converted FP64 runs, and \textbf{9.00\%} of FP32 
runs, highlighting differences in numerical consistency 
across precision settings and compiler environments. 

\begin{mdframed}[backgroundcolor=light-gray] 
\textbf{Answer to Q1:}
Our approach is effective at identifying compiler-induced numerical 
inconsistencies between NVIDIA and AMD GPUs. As shown in 
Table~\ref{tab:sum_exp_res}, we observed discrepancies in 0.98\% of the 
247,500 runs conducted under FP64 settings. In tests using 
\HIPIFY-converted FP64 settings, 
discrepancies increased to 1.10\%, indicating a significant difference. 
For FP32 tests, discrepancies were found in 9\% of the 157,600 runs. 
\textbf{These results  
demonstrate the effectiveness of our approach in detecting numerical inconsistencies.}
\end{mdframed}

\begin{table}
\centering
\caption{Summary of Experimental Results}
\label{tab:sum_exp_res}
\begin{adjustbox}{max width=\columnwidth}
\begin{tabular}{lccc}
\toprule
\toprule
\scriptsize
\textbf{Metric} & \textbf{FP64} & \textbf{FP64 with \hipify} & \textbf{FP32} \\
\midrule
Total Programs & 3,540 & 3,540 & 2,840 \\
Total Runs per Option per Compiler & 24,750 & 24,750 & 15,760 \\
Total Runs per Option & 49,500 & 49,500 & 31,520 \\
\textbf{Total Runs} & \textbf{247,500} & \textbf{247,500} & \textbf{157,600} \\
\midrule
Runs on NVCC & 123,750 & 123,750 & 78,800 \\
Runs on HIPCC & 123,750 & 123,750 & 78,800 \\
\midrule
\textbf{Total Discrepancies} & \textbf{2,426} & \textbf{2,716} & \textbf{14,188} \\
\textit{Total Discrepancies (\% of Total Runs)} & \textit{0.98\%} & \textit{1.10\%} & \textit{9.00\%} \\
\bottomrule
\bottomrule
\end{tabular}
\end{adjustbox}
\end{table}

\subsubsection{FP64 Tests}
Table~\ref{tab:f64_res_opt} categorizes discrepancies 
for FP64 tests by optimization level. 
While \texttt{O3} with \texttt{-ffast-math} (\texttt{O3\_FM}) exhibited the 
highest discrepancies at \texttt{519}, the \texttt{O0} setting 
also had a significant count at \texttt{440}, suggesting 
that both aggressive and no optimizations can lead to variations. 
The Number vs. Number discrepancies were the most frequent 
across all optimization settings, indicating persistent 
challenges in achieving numerical consistency. 
Table~\ref{tab:f64_res_adj} presents adjacency matrices 
for various optimization levels.

\begin{table*}
\centering
\caption{Discrepancies per Optimization Option for FP64 tests}
\label{tab:f64_res_opt}
\begin{adjustbox}{width=0.9\textwidth}
\begin{tabular}{l@{}p{1.5cm}@{}p{1.5cm}@{}p{1.5cm}@{}p{1.5cm}@{}p{1.5cm}@{}p{1.5cm}@{}p{1.5cm}@{}p{1.5cm}@{}}
\toprule
\textbf{Opt Flags } & \textbf{Disc. Count } & \textbf{NaN, Inf } & \textbf{NaN, Zero} & \textbf{NaN, Num} & \textbf{Inf, Zero} & \textbf{Inf, Num} & \textbf{Num, Zero} & \textbf{Num, Num} \\
\midrule
\textbf{O0} & \raisebox{0.5mm}{\colorbox{cyan!30}{\phantom{0}}} \textbf{440} & \raisebox{0.5mm}{\colorbox{cyan!10}{\phantom{0}}} 7 & \raisebox{0.5mm}{\colorbox{cyan!0}{\framebox(4,4){}}} 0 & \raisebox{0.5mm}{\colorbox{cyan!0}{\framebox(4,4){}}} 0 & \raisebox{0.5mm}{\colorbox{cyan!15}{\phantom{0}}} 24 & \raisebox{0.5mm}{\colorbox{cyan!20}{\phantom{0}}} 46 & \raisebox{0.5mm}{\colorbox{cyan!20}{\phantom{0}}} 10 & \raisebox{0.5mm}{\colorbox{cyan!70}{\phantom{0}}} \textbf{353} \\
\textbf{O1} & \raisebox{0.5mm}{\colorbox{cyan!33}{\phantom{0}}} \textbf{489} & \raisebox{0.5mm}{\colorbox{cyan!13}{\phantom{0}}} 22 & \raisebox{0.5mm}{\colorbox{cyan!0}{\framebox(4,4){}}} 0 & \raisebox{0.5mm}{\colorbox{cyan!0}{\framebox(4,4){}}} 0 & \raisebox{0.5mm}{\colorbox{cyan!15}{\phantom{0}}} 24 & \raisebox{0.5mm}{\colorbox{cyan!20}{\phantom{0}}} 46 & \raisebox{0.5mm}{\colorbox{cyan!20}{\phantom{0}}} 10 & \raisebox{0.5mm}{\colorbox{cyan!67}{\phantom{0}}} \textbf{387} \\
\textbf{O2} & \raisebox{0.5mm}{\colorbox{cyan!33}{\phantom{0}}} \textbf{489} & \raisebox{0.5mm}{\colorbox{cyan!13}{\phantom{0}}} 22 & \raisebox{0.5mm}{\colorbox{cyan!0}{\framebox(4,4){}}} 0 & \raisebox{0.5mm}{\colorbox{cyan!0}{\framebox(4,4){}}} 0 & \raisebox{0.5mm}{\colorbox{cyan!15}{\phantom{0}}} 24 & \raisebox{0.5mm}{\colorbox{cyan!20}{\phantom{0}}} 46 & \raisebox{0.5mm}{\colorbox{cyan!20}{\phantom{0}}} 10 & \raisebox{0.5mm}{\colorbox{cyan!67}{\phantom{0}}} \textbf{387} \\
\textbf{O3} & \raisebox{0.5mm}{\colorbox{cyan!33}{\phantom{0}}} \textbf{489} & \raisebox{0.5mm}{\colorbox{cyan!13}{\phantom{0}}} 22 & \raisebox{0.5mm}{\colorbox{cyan!0}{\framebox(4,4){}}} 0 & \raisebox{0.5mm}{\colorbox{cyan!0}{\framebox(4,4){}}} 0 & \raisebox{0.5mm}{\colorbox{cyan!15}{\phantom{0}}} 24 & \raisebox{0.5mm}{\colorbox{cyan!20}{\phantom{0}}} 46 & \raisebox{0.5mm}{\colorbox{cyan!20}{\phantom{0}}} 10 & \raisebox{0.5mm}{\colorbox{cyan!67}{\phantom{0}}} \textbf{387} \\
\textbf{O3\_FM} & \raisebox{0.5mm}{\colorbox{cyan!35}{\phantom{0}}} \textbf{519} & \raisebox{0.5mm}{\colorbox{cyan!23}{\phantom{0}}} \textbf{32} & \raisebox{0.5mm}{\colorbox{cyan!0}{\framebox(4,4){}}} \textbf{0} & \raisebox{0.5mm}{\colorbox{cyan!0}{\framebox(4,4){}}} \textbf{0} & \raisebox{0.5mm}{\colorbox{cyan!15}{\phantom{0}}} \textbf{24} & \raisebox{0.5mm}{\colorbox{cyan!23}{\phantom{0}}} \textbf{52} & \raisebox{0.5mm}{\colorbox{cyan!20}{\phantom{0}}} \textbf{10} & \raisebox{0.5mm}{\colorbox{cyan!68}{\phantom{0}}} \textbf{401} \\
\textbf{Total} & \raisebox{0.5mm}{\colorbox{cyan!100}{\phantom{0}}} \textbf{2,426} & \raisebox{0.5mm}{\colorbox{cyan!40}{\phantom{0}}} \textbf{105} & \raisebox{0.5mm}{\colorbox{cyan!0}{\framebox(4,4){}}} \textbf{0} & \raisebox{0.5mm}{\colorbox{cyan!0}{\framebox(4,4){}}} \textbf{0} & \raisebox{0.5mm}{\colorbox{cyan!25}{\phantom{0}}} \textbf{120} & \raisebox{0.5mm}{\colorbox{cyan!39}{\phantom{0}}} \textbf{236} & \raisebox{0.5mm}{\colorbox{cyan!37}{\phantom{0}}} \textbf{50} & \raisebox{0.5mm}{\colorbox{cyan!86}{\phantom{0}}} \textbf{1,915} \\
\bottomrule
\end{tabular}
\end{adjustbox}

\vspace{0.5em} 

\begin{adjustbox}{max width=\textwidth}
\begin{tabular}{ccccccccccc}
\toprule
0\% & 10\% & 20\% & 30\% & 40\% & 50\% & 60\% & 70\% & 80\% & 90\% & 100\% \\
\cellcolor{cyan!0}\phantom{0} & \cellcolor{cyan!10}\phantom{0} & \cellcolor{cyan!20}\phantom{0} & \cellcolor{cyan!30}\phantom{0} & \cellcolor{cyan!40}\phantom{0} & \cellcolor{cyan!50}\phantom{0} & \cellcolor{cyan!60}\phantom{0} & \cellcolor{cyan!70}\phantom{0} & \cellcolor{cyan!80}\phantom{0} & \cellcolor{cyan!90}\phantom{0} & \cellcolor{cyan!100}\phantom{0} \\
\bottomrule
\end{tabular}
\end{adjustbox}
\end{table*}

\begin{table}
\centering
\caption{Adjacency Matrices for Different Optimization Levels for FP64 tests}
\label{tab:f64_res_adj}
\begin{adjustbox}{max width=\columnwidth}
\begin{tabular}{lccccc}
\toprule
\toprule
\textbf{Opt Flags} & \textbf{NVCC\textbackslash HIPCC} & \textbf{(±) NaN} & \textbf{(±) Inf} & \textbf{(±) Zero} & \textbf{Num} \\
\midrule
\midrule
\multirow{4}{*}{\textbf{O0}} & \textbf{(±) NaN} & — & 7, 0 & 0, 0 & 0, 0 \\
& \textbf{(±) Inf} & — & — & 8, 16 & 4, 42 \\
& \textbf{(±) Zero} & — & — & — & 0, 10 \\
& \textbf{Num} & — & — & — & \textbf{353, 353} \\
\midrule
\midrule
\multirow{4}{*}{\textbf{O1}} & \textbf{(±) NaN} & — & 19, 3 & 0, 0 & 0, 0 \\
& \textbf{(±) Inf} & — & — & 8, 16 & 4, 42 \\
& \textbf{(±) Zero} & — & — & — & 0, 10 \\
& \textbf{Num} & — & — & — & \textbf{387, 387} \\
\midrule
\midrule
\multirow{4}{*}{\textbf{O2}} & \textbf{(±) NaN} & — & 19, 3 & 0, 0 & 0, 0 \\
& \textbf{(±) Inf} & — & — & 8, 16 & 4, 42 \\
& \textbf{(±) Zero} & — & — & — & 0, 10 \\
& \textbf{Num} & — & — & — & \textbf{387, 387} \\
\midrule
\midrule
\multirow{4}{*}{\textbf{O3}} & \textbf{(±) NaN} & — & 19, 3 & 0, 0 & 0, 0 \\
& \textbf{(±) Inf} & — & — & 8, 16 & 4, 42 \\
& \textbf{(±) Zero} & — & — & — & 0, 10 \\
& \textbf{Num} & — & — & — & \textbf{387, 387} \\
\midrule
\midrule
\multirow{4}{*}{\textbf{O3\_FM}} & \textbf{(±) NaN} & — & \textbf{29, 3} & \textbf{0, 0} & \textbf{0, 0} \\
& \textbf{(±) Inf} & — & — & \textbf{8, 16} & \textbf{4, 48} \\
& \textbf{(±) Zero} & — & — & — & \textbf{0, 10} \\
& \textbf{Num} & — & — & — & \textbf{401, 401} \\
\bottomrule
\bottomrule
\end{tabular}
\end{adjustbox}
\end{table}

\subsubsection{\HIPIFY-Converted FP64 Tests}
This experiment allowed us to test whether 
the \HIPIFY tool introduces additional inconsistencies 
during the conversion from CUDA to HIP. We found that 
\HIPIFY introduces discrepancies 
in floating-point computations between NVIDIA and AMD architectures. 
Table~\ref{tab:hip_res_opt} categorizes the discrepancies observed for \HIPIFY-converted 
tests across different optimization levels. 
The results reveal that \texttt{O3} with \texttt{-ffast-math} (\texttt{O3\_FM}) 
exhibits the highest discrepancy count at \texttt{575}, 
while \texttt{O0} records \texttt{494} discrepancies. 
Discrepancies in the Number vs. Number category remain 
predominant across all optimization levels, 
like the previous \texttt{FP64} tests. Table~\ref{tab:hip_res_adj} 
presents adjacency matrices for various optimization levels.

\begin{table*}
\centering
\caption{Discrepancies per Optimization Option for \HIPIFY converted FP64}
\label{tab:hip_res_opt}
\begin{adjustbox}{width=0.9\textwidth}
\begin{tabular}{l@{}p{1.5cm}@{}p{1.5cm}@{}p{1.5cm}@{}p{1.5cm}@{}p{1.5cm}@{}p{1.5cm}@{}p{1.5cm}@{}p{1.5cm}@{}}
\toprule
\textbf{Opt Flags} & \textbf{Disc. Count} & \textbf{NaN, Inf} & \textbf{NaN, Zero} & \textbf{NaN, Num} & \textbf{Inf, Zero} & \textbf{Inf, Num} & \textbf{Num, Zero} & \textbf{Num, Num} \\
\midrule
\textbf{O0} & \raisebox{0.5mm}{\colorbox{cyan!18}{\phantom{0}}} \textbf{494} & \raisebox{0.5mm}{\colorbox{cyan!1}{\phantom{0}}} 3 & \raisebox{0.5mm}{\colorbox{cyan!0}{\framebox(4,4){}}} 0 & \raisebox{0.5mm}{\colorbox{cyan!0}{\framebox(4,4){}}} 0 & \raisebox{0.5mm}{\colorbox{cyan!4}{\phantom{0}}} 12 & \raisebox{0.5mm}{\colorbox{cyan!7}{\phantom{0}}} 20 & \raisebox{0.5mm}{\colorbox{cyan!7}{\phantom{0}}} 20 & \raisebox{0.5mm}{\colorbox{cyan!16}{\phantom{0}}} \textbf{439} \\
\textbf{O1} & \raisebox{0.5mm}{\colorbox{cyan!20}{\phantom{0}}} \textbf{549} & \raisebox{0.5mm}{\colorbox{cyan!8}{\phantom{0}}} 23 & \raisebox{0.5mm}{\colorbox{cyan!0}{\framebox(4,4){}}} 0 & \raisebox{0.5mm}{\colorbox{cyan!0}{\framebox(4,4){}}} 0 & \raisebox{0.5mm}{\colorbox{cyan!4}{\phantom{0}}} 12 & \raisebox{0.5mm}{\colorbox{cyan!7}{\phantom{0}}} 20 & \raisebox{0.5mm}{\colorbox{cyan!7}{\phantom{0}}} 20 & \raisebox{0.5mm}{\colorbox{cyan!17}{\phantom{0}}} \textbf{474} \\
\textbf{O2} & \raisebox{0.5mm}{\colorbox{cyan!20}{\phantom{0}}} \textbf{549} & \raisebox{0.5mm}{\colorbox{cyan!8}{\phantom{0}}} 23 & \raisebox{0.5mm}{\colorbox{cyan!0}{\framebox(4,4){}}} 0 & \raisebox{0.5mm}{\colorbox{cyan!0}{\framebox(4,4){}}} 0 & \raisebox{0.5mm}{\colorbox{cyan!4}{\phantom{0}}} 12 & \raisebox{0.5mm}{\colorbox{cyan!7}{\phantom{0}}} 20 & \raisebox{0.5mm}{\colorbox{cyan!7}{\phantom{0}}} 20 & \raisebox{0.5mm}{\colorbox{cyan!17}{\phantom{0}}} \textbf{474} \\
\textbf{O3} & \raisebox{0.5mm}{\colorbox{cyan!20}{\phantom{0}}} \textbf{549} & \raisebox{0.5mm}{\colorbox{cyan!8}{\phantom{0}}} 23 & \raisebox{0.5mm}{\colorbox{cyan!0}{\framebox(4,4){}}} 0 & \raisebox{0.5mm}{\colorbox{cyan!0}{\framebox(4,4){}}} 0 & \raisebox{0.5mm}{\colorbox{cyan!4}{\phantom{0}}} 12 & \raisebox{0.5mm}{\colorbox{cyan!7}{\phantom{0}}} 20 & \raisebox{0.5mm}{\colorbox{cyan!7}{\phantom{0}}} 20 & \raisebox{0.5mm}{\colorbox{cyan!17}{\phantom{0}}} \textbf{474} \\
\textbf{O3\_FM} & \raisebox{0.5mm}{\colorbox{cyan!21}{\phantom{0}}} \textbf{575} & \raisebox{0.5mm}{\colorbox{cyan!13}{\phantom{0}}} \textbf{36} & \raisebox{0.5mm}{\colorbox{cyan!0}{\framebox(4,4){}}} \textbf{0} & \raisebox{0.5mm}{\colorbox{cyan!0}{\framebox(4,4){}}} \textbf{0} & \raisebox{0.5mm}{\colorbox{cyan!4}{\phantom{0}}} \textbf{12} & \raisebox{0.5mm}{\colorbox{cyan!10}{\phantom{0}}} \textbf{27} & \raisebox{0.5mm}{\colorbox{cyan!7}{\phantom{0}}} \textbf{20} & \raisebox{0.5mm}{\colorbox{cyan!18}{\phantom{0}}} \textbf{480} \\
\textbf{Total} & \raisebox{0.5mm}{\colorbox{cyan!100}{\phantom{0}}} \textbf{2,716} & \raisebox{0.5mm}{\colorbox{cyan!40}{\phantom{0}}} \textbf{108} & \raisebox{0.5mm}{\colorbox{cyan!0}{\framebox(4,4){}}} \textbf{0} & \raisebox{0.5mm}{\colorbox{cyan!0}{\framebox(4,4){}}} \textbf{0} & \raisebox{0.5mm}{\colorbox{cyan!2}{\phantom{0}}} \textbf{60} & \raisebox{0.5mm}{\colorbox{cyan!4}{\phantom{0}}} \textbf{107} & \raisebox{0.5mm}{\colorbox{cyan!3}{\phantom{0}}} \textbf{100} & \raisebox{0.5mm}{\colorbox{cyan!86}{\phantom{0}}} \textbf{2,341} \\
\bottomrule
\end{tabular}
\end{adjustbox}

\vspace{0.5em}

\begin{adjustbox}{max width=\textwidth}
\begin{tabular}{ccccccccccc}
\toprule
0\% & 10\% & 20\% & 30\% & 40\% & 50\% & 60\% & 70\% & 80\% & 90\% & 100\% \\
\cellcolor{cyan!0}\phantom{0} & \cellcolor{cyan!10}\phantom{0} & \cellcolor{cyan!20}\phantom{0} & \cellcolor{cyan!30}\phantom{0} & \cellcolor{cyan!40}\phantom{0} & \cellcolor{cyan!50}\phantom{0} & \cellcolor{cyan!60}\phantom{0} & \cellcolor{cyan!70}\phantom{0} & \cellcolor{cyan!80}\phantom{0} & \cellcolor{cyan!90}\phantom{0} & \cellcolor{cyan!100}\phantom{0} \\
\bottomrule
\end{tabular}
\end{adjustbox}
\end{table*}

\begin{table}
\centering
\caption{Adjacency Matrices for Different Optimization Levels for \HIPIFY converted FP64}
\label{tab:hip_res_adj}
\begin{adjustbox}{max width=\columnwidth}
\begin{tabular}{lccccc}
\toprule
\toprule
\textbf{Opt Flags} & \textbf{NVCC\textbackslash HIPCC} & \textbf{(±) NaN} & \textbf{(±) Inf} & \textbf{(±) Zero} & \textbf{Num} \\
\midrule
\midrule
\multirow{4}{*}{\textbf{O0}} & \textbf{(±) NaN} & — & 0, 3 & 0, 0 & 0, 0 \\
& \textbf{(±) Inf} & — & — & 8, 4 & 2, 18 \\
& \textbf{(±) Zero} & — & — & — & 0, 20 \\
& \textbf{Num} & — & — & — & \textbf{439, 439} \\
\midrule
\midrule
\multirow{4}{*}{\textbf{O1}} & \textbf{(±) NaN} & — & 17, 6 & 0, 0 & 0, 0 \\
& \textbf{(±) Inf} & — & — & 8, 4 & 2, 18 \\
& \textbf{(±) Zero} & — & — & — & 0, 20 \\
& \textbf{Num} & — & — & — & \textbf{474, 474} \\
\midrule
\midrule
\multirow{4}{*}{\textbf{O2}} & \textbf{(±) NaN} & — & 17, 6 & 0, 0 & 0, 0 \\
& \textbf{(±) Inf} & — & — & 8, 4 & 2, 18 \\
& \textbf{(±) Zero} & — & — & — & 0, 20 \\
& \textbf{Num} & — & — & — & \textbf{474, 474} \\
\midrule
\midrule
\multirow{4}{*}{\textbf{O3}} & \textbf{(±) NaN} & — & 17, 6 & 0, 0 & 0, 0 \\
& \textbf{(±) Inf} & — & — & 8, 4 & 2, 18 \\
& \textbf{(±) Zero} & — & — & — & 0, 20 \\
& \textbf{Num} & — & — & — & \textbf{474, 474} \\
\midrule
\midrule
\multirow{4}{*}{\textbf{O3\_FM}} & \textbf{(±) NaN} & — & \textbf{30, 6} & \textbf{0, 0} & \textbf{0, 0} \\
& \textbf{(±) Inf} & — & — & \textbf{8, 4} & \textbf{2, 25} \\
& \textbf{(±) Zero} & — & — & — & \textbf{0, 20} \\
& \textbf{Num} & — & — & — & \textbf{480, 480} \\
\bottomrule
\bottomrule
\end{tabular}
\end{adjustbox}
\end{table}

\subsubsection{FP32 Tests}
Single-precision floating-point computations are prevalent 
in machine learning applications running on 32-bit systems. 
Table~\ref{tab:f32_res_opt} categorizes the discrepancies 
observed per optimization flag for FP32 tests, 
revealing a stark contrast in numerical consistency compared to FP64 tests. 
The results show that the \texttt{O3} with \texttt{-ffast-math} (\texttt{O3\_FM}) option yielded the highest discrepancy count at \texttt{13,877}, compared to only \texttt{45} 
discrepancies in the \texttt{O0} setting. This significant increase 
highlights the aggressive nature of \texttt{O3\_FM} optimizations, 
which prioritize speed and may sacrifice numerical precision. 
Notably, the Number vs. Number category accounted for the majority of discrepancies, 
indicating persistent challenges in maintaining consistency. 
Table~\ref{tab:f32_res_adj} provides a comprehensive analysis 
of discrepancies across various optimization levels.

\begin{table*}
\centering
\caption{Discrepancies per Optimization Option for FP32 tests}
\label{tab:f32_res_opt}
\begin{adjustbox}{width=0.9\textwidth}
\begin{tabular}{l@{}p{1.5cm}@{}p{1.5cm}@{}p{1.5cm}@{}p{1.5cm}@{}p{1.5cm}@{}p{1.5cm}@{}p{1.5cm}@{}p{1.5cm}@{}}
\toprule
\textbf{Opt Flags} & \textbf{Disc. Count} & \textbf{NaN, Inf} & \textbf{NaN, Zero} & \textbf{NaN, Num} & \textbf{Inf, Zero} & \textbf{Inf, Num} & \textbf{Num, Zero} & \textbf{Num, Num} \\
\midrule
\textbf{O0} & \raisebox{0.5mm}{\colorbox{cyan!1}{\phantom{0}}} \textbf{45} & \raisebox{0.5mm}{\colorbox{cyan!3}{\phantom{0}}} 5 & \raisebox{0.5mm}{\colorbox{cyan!0}{\framebox(4,4){}}} 0 & \raisebox{0.5mm}{\colorbox{cyan!0}{\framebox(4,4){}}} 0 & \raisebox{0.5mm}{\colorbox{cyan!0}{\framebox(4,4){}}} 0 & \raisebox{0.5mm}{\colorbox{cyan!0}{\framebox(4,4){}}} 0 & \raisebox{0.5mm}{\colorbox{cyan!0}{\framebox(4,4){}}} 0 & \raisebox{0.5mm}{\colorbox{cyan!1}{\phantom{0}}} \textbf{40} \\
\textbf{O1} & \raisebox{0.5mm}{\colorbox{cyan!6}{\phantom{0}}} \textbf{86} & \raisebox{0.5mm}{\colorbox{cyan!3}{\phantom{0}}} 24 & \raisebox{0.5mm}{\colorbox{cyan!0}{\framebox(4,4){}}} 0 & \raisebox{0.5mm}{\colorbox{cyan!0}{\framebox(4,4){}}} 6 & \raisebox{0.5mm}{\colorbox{cyan!0}{\framebox(4,4){}}} 0 & \raisebox{0.5mm}{\colorbox{cyan!0}{\framebox(4,4){}}} 0 & \raisebox{0.5mm}{\colorbox{cyan!0}{\framebox(4,4){}}} 0 & \raisebox{0.5mm}{\colorbox{cyan!6}{\phantom{0}}} \textbf{56} \\
\textbf{O2} & \raisebox{0.5mm}{\colorbox{cyan!6}{\phantom{0}}} \textbf{90} & \raisebox{0.5mm}{\colorbox{cyan!4}{\phantom{0}}} 28 & \raisebox{0.5mm}{\colorbox{cyan!0}{\framebox(4,4){}}} 0 & \raisebox{0.5mm}{\colorbox{cyan!0}{\framebox(4,4){}}} 6 & \raisebox{0.5mm}{\colorbox{cyan!0}{\framebox(4,4){}}} 0 & \raisebox{0.5mm}{\colorbox{cyan!0}{\framebox(4,4){}}} 0 & \raisebox{0.5mm}{\colorbox{cyan!0}{\framebox(4,4){}}} 0 & \raisebox{0.5mm}{\colorbox{cyan!6}{\phantom{0}}} \textbf{56} \\
\textbf{O3} & \raisebox{0.5mm}{\colorbox{cyan!6}{\phantom{0}}} \textbf{90} & \raisebox{0.5mm}{\colorbox{cyan!4}{\phantom{0}}} 28 & \raisebox{0.5mm}{\colorbox{cyan!0}{\framebox(4,4){}}} 0 & \raisebox{0.5mm}{\colorbox{cyan!0}{\framebox(4,4){}}} 6 & \raisebox{0.5mm}{\colorbox{cyan!0}{\framebox(4,4){}}} 0 & \raisebox{0.5mm}{\colorbox{cyan!0}{\framebox(4,4){}}} 0 & \raisebox{0.5mm}{\colorbox{cyan!0}{\framebox(4,4){}}} 0 & \raisebox{0.5mm}{\colorbox{cyan!6}{\phantom{0}}} \textbf{56} \\
\textbf{O3\_FM} & \raisebox{0.5mm}{\colorbox{cyan!98}{\phantom{0}}} \textbf{13877} & \raisebox{0.5mm}{\colorbox{cyan!16}{\phantom{0}}} \textbf{2328} & \raisebox{0.5mm}{\colorbox{cyan!4}{\phantom{0}}} \textbf{619} & \raisebox{0.5mm}{\colorbox{cyan!15}{\phantom{0}}} \textbf{2134} & \raisebox{0.5mm}{\colorbox{cyan!1}{\phantom{0}}} \textbf{259} & \raisebox{0.5mm}{\colorbox{cyan!7}{\phantom{0}}} \textbf{1056} & \raisebox{0.5mm}{\colorbox{cyan!38}{\phantom{0}}} \textbf{5369} & \raisebox{0.5mm}{\colorbox{cyan!15}{\phantom{0}}} \textbf{2112} \\
\textbf{Total} & \raisebox{0.5mm}{\colorbox{cyan!100}{\phantom{0}}} \textbf{14188} & \raisebox{0.5mm}{\colorbox{cyan!17}{\phantom{0}}} \textbf{2413} & \raisebox{0.5mm}{\colorbox{cyan!4}{\phantom{0}}} \textbf{619} & \raisebox{0.5mm}{\colorbox{cyan!15}{\phantom{0}}} \textbf{2152} & \raisebox{0.5mm}{\colorbox{cyan!1}{\phantom{0}}} \textbf{259} & \raisebox{0.5mm}{\colorbox{cyan!7}{\phantom{0}}} \textbf{1056} & \raisebox{0.5mm}{\colorbox{cyan!38}{\phantom{0}}} \textbf{5369} & \raisebox{0.5mm}{\colorbox{cyan!16}{\phantom{0}}} \textbf{2320} \\
\bottomrule
\end{tabular}
\end{adjustbox}

\vspace{0.5em}

\begin{adjustbox}{max width=\textwidth}
\begin{tabular}{ccccccccccc}
\toprule
0\% & 10\% & 20\% & 30\% & 40\% & 50\% & 60\% & 70\% & 80\% & 90\% & 100\% \\
\cellcolor{cyan!0}\phantom{0} & \cellcolor{cyan!10}\phantom{0} & \cellcolor{cyan!20}\phantom{0} & \cellcolor{cyan!30}\phantom{0} & \cellcolor{cyan!40}\phantom{0} & \cellcolor{cyan!50}\phantom{0} & \cellcolor{cyan!60}\phantom{0} & \cellcolor{cyan!70}\phantom{0} & \cellcolor{cyan!80}\phantom{0} & \cellcolor{cyan!90}\phantom{0} & \cellcolor{cyan!100}\phantom{0} \\
\bottomrule
\end{tabular}
\end{adjustbox}
\end{table*}

\begin{table}
\centering
\caption{Adjacency Matrices for Different Optimization Levels for FP32 tests}
\label{tab:f32_res_adj}
\begin{adjustbox}{max width=\columnwidth}
\begin{tabular}{lccccc}
\toprule
\toprule
\textbf{Opt Flags} & \textbf{NVCC\textbackslash HIPCC} & \textbf{(±) NaN} & \textbf{(±) Inf} & \textbf{(±) Zero} & \textbf{Num} \\
\midrule
\midrule
\multirow{4}{*}{\textbf{O0}} & \textbf{(±) NaN} & — & 0, 5 & 0, 0 & 0, 0 \\
& \textbf{(±) Inf} & — & — & 0, 0 & 0, 0 \\
& \textbf{(±) Zero} & — & — & — & 0, 0 \\
& \textbf{Num} & — & — & — & \textbf{40, 40} \\
\midrule
\midrule
\multirow{4}{*}{\textbf{O1}} & \textbf{(±) NaN} & — & 19, 5 & 0, 0 & 6, 0 \\
& \textbf{(±) Inf} & — & — & 0, 0 & 0, 0 \\
& \textbf{(±) Zero} & — & — & — & 0, 0 \\
& \textbf{Num} & — & — & — & \textbf{56, 56} \\
\midrule
\midrule
\multirow{4}{*}{\textbf{O2}} & \textbf{(±) NaN} & — & 23, 5 & 0, 0 & 6, 0 \\
& \textbf{(±) Inf} & — & — & 0, 0 & 0, 0 \\
& \textbf{(±) Zero} & — & — & — & 0, 0 \\
& \textbf{Num} & — & — & — & \textbf{56, 56} \\
\midrule
\midrule
\multirow{4}{*}{\textbf{O3}} & \textbf{(±) NaN} & — & 23, 5 & 0, 0 & 6, 0 \\
& \textbf{(±) Inf} & — & — & 0, 0 & 0, 0 \\
& \textbf{(±) Zero} & — & — & — & 0, 0 \\
& \textbf{Num} & — & — & — & \textbf{56, 56} \\
\midrule
\midrule
\multirow{4}{*}{\textbf{O3\_fast\_math}} & \textbf{(±) NaN} & — & \textbf{101, 2227} & \textbf{113, 506} & \textbf{139, 1995} \\
& \textbf{(±) Inf} & — & — & \textbf{190, 69} & \textbf{197, 859} \\
& \textbf{(±) Zero} & — & — & — & \textbf{161, 5208} \\
& \textbf{Num} & — & — & — & \textbf{2112, 2112} \\
\bottomrule
\bottomrule
\end{tabular}
\end{adjustbox}
\end{table}

\begin{mdframed}[backgroundcolor=light-gray] 
\textbf{Answer to Q2:} We identified seven classes of numerical 
inconsistencies.
These discrepancies arise from differences in floating-point computations 
between NVIDIA and AMD GPUs. In our FP64 tests, we did not 
observe NaN vs. Zero and NaN vs. Number discrepancies, indicating that 
a larger search space may be necessary to detect them, if they exist. 
In FP32 tests, discrepancies such as NaN vs. Inf, Nan vs. Num, and Num vs. Zero
predominantly occurred under the \texttt{-O3} flag with \texttt{-ffast-math} enabled. 
\textbf{Overall, we observed all classes of inconsistencies across various tests}.
\end{mdframed}


\subsection{Case Studies}
\label{sec:cases}
in this section, we focus on three interesting cases.

\subsubsection{Case Study 1: Inconsistency in Real-Valued Results}
\begin{figure}
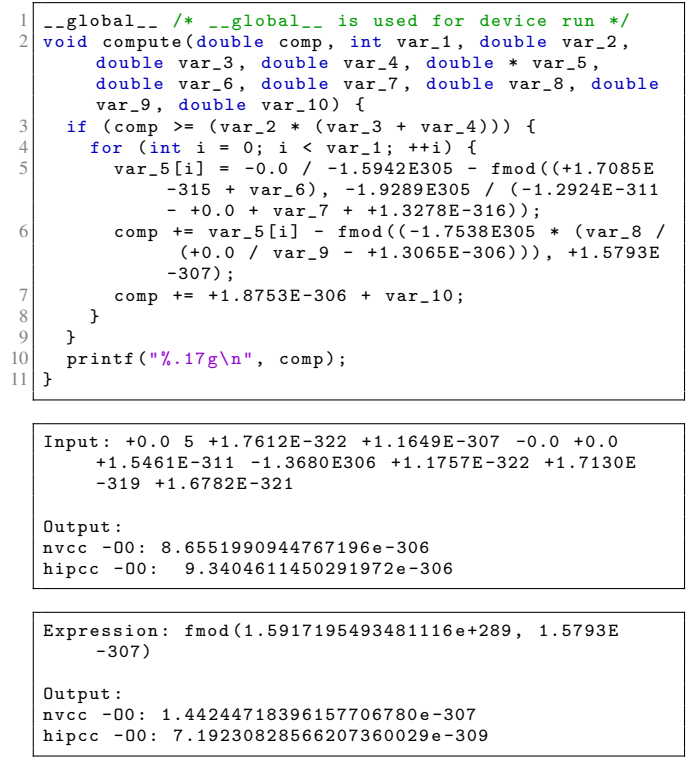

\centering
\begin{lstlisting}[language=C]
__global__ /* __global__ is used for device run */
void compute(double comp, int var_1, double var_2, double var_3, double var_4, double * var_5, double var_6, double var_7, double var_8, double var_9, double var_10) {
  if (comp >= (var_2 * (var_3 + var_4))) {
    for (int i = 0; i < var_1; ++i) {
      var_5[i] = -0.0 / -1.5942E305 - fmod((+1.7085E-315 + var_6), -1.9289E305 / (-1.2924E-311 - +0.0 + var_7 + +1.3278E-316));
      comp += var_5[i] - fmod((-1.7538E305 * (var_8 / (+0.0 / var_9 - +1.3065E-306))), +1.5793E-307);
      comp += +1.8753E-306 + var_10;
    }
  }
  printf("%.17g\n", comp);
}
\end{lstlisting}
\begin{lstlisting}[language=C, numbers=none]
Input: +0.0 5 +1.7612E-322 +1.1649E-307 -0.0 +0.0 +1.5461E-311 -1.3680E306 +1.1757E-322 +1.7130E-319 +1.6782E-321

Output:
nvcc -O0:	8.6551990944767196e-306
hipcc -O0:	9.3404611450291972e-306
\end{lstlisting}
\begin{lstlisting}[language=C, numbers=none]
Expression: fmod(1.5917195493481116e+289, 1.5793E-307)

Output:
nvcc -O0: 1.44244718396157706780e-307
hipcc -O0: 7.19230828566207360029e-309
\end{lstlisting}
\caption{Small Numerical Variation with No Optimization (\texttt{-O0})}
\label{fig:case_study_1}
\end{figure}

During our experiments, \varity generated the code shown in 
Figure~\ref{fig:case_study_1} and it comprises three parts. 
The first part shows a code snippet. Upon inspection, 
the \textbf{compute} method performs several floating-point operations, 
including division and the use of the \texttt{fmod} function for remainders. 
A conditional statement determines whether the code enters a loop, 
where we identified variations in the results generated within the loop.

In the second part, we see the failure-inducing input. 
For the same input,
the \texttt{nvcc} compiler produced an output of \texttt{8.6551990944767196e-306}, 
while the \texttt{hipcc} compiler produced \texttt{9.3404611450291972e-306}. 
To determine the cause of this discrepancy, we analyzed the 
intermediate results and the assembly code generated by the compilers. 

We found that until the condition was satisfied and the loop started, 
there were no issues with this input. During the first iteration, 
the expression \texttt{(-1.7538E305 * (var\_8 / (+0.0 / var\_9 - +1.3065E-306)))} 
was computed, resulting in the same value on both devices: \texttt{1.5917195493481116e+289}. 
The discrepancy arose when this result was passed to the \texttt{fmod} function. 
The \texttt{nvcc} compiler computed \texttt{fmod(1.5917195493481116e+289, 1.5793E-307)} 
as \texttt{1.4424471839615771e-307}, whereas the \texttt{hipcc} compiler 
produced \texttt{7.1923082856620736e-309}. This small numerical 
difference then propagated through subsequent iterations, 
magnified with each loop iteration, and by the end of the loop, 
these differences compounded, resulting in significantly 
different final outputs on the two devices, as shown in the third part of the Figure.

Further examination of the assembly code revealed that the \texttt{fmod} 
function is implemented differently on the two platforms. 
For the \texttt{hipcc} compiler, the \texttt{fmod} function is 
a clang-style function named \texttt{\_\_ocml\_fmod\_f64}, which is 
called each time the function is invoked. This implementation is part 
of the AMD GPU Instruction Set Architecture (ISA). In contrast, 
the \texttt{nvcc} compiler uses a combination of floating-point 
arithmetic and bitwise manipulation within its SASS (String Assembler) and PTX (Parallel Thread Execution) assembly 
languages to implement the \texttt{fmod} function. 

Interestingly, out of ten randomly generated inputs, 
only this specific input created a discrepancy. 
Other inputs, such as \texttt{"-0.0 5 +0.0 +1.2150E-306 +1.2318E224 +1.8418E306 +1.6483E-306 -1.2836E214 -1.0053E305 +1.3789E213 +0.0"} or \texttt{"-1.2640E305 5 -0.0 -1.3396E-322 +1.7693E-307 +1.3050E305 +1.4789E-316 -1.1350E-165 +1.7826E305 +0.0 -1.0755E-317"}, 
produced consistent results between the two platforms. 
This case study highlights the challenges in achieving 
consistent floating-point computations across two prominent 
GPU architectures due to differences in \textit{math functions}.

\subsubsection{Case Study 2: Inconsistencies in Infinity-Valued Results}
\begin{figure}
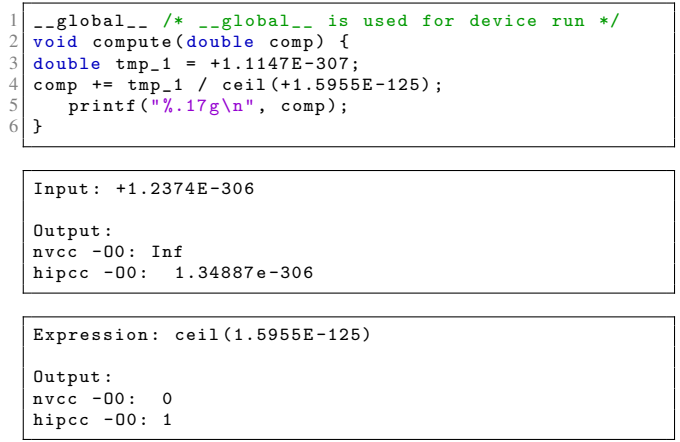

\centering
\begin{lstlisting}[language=C]
__global__ /* __global__ is used for device run */
void compute(double comp) {
double tmp_1 = +1.1147E-307;
comp += tmp_1 / ceil(+1.5955E-125);
   printf("%.17g\n", comp);
}
\end{lstlisting}
\begin{lstlisting}[language=C, numbers=none]
Input: +1.2374E-306

Output:
nvcc -O0:	Inf
hipcc -O0:	1.34887e-306
\end{lstlisting}
\begin{lstlisting}[language=C, numbers=none]
Expression: ceil(1.5955E-125)

Output:
nvcc -O0:  0
hipcc -O0: 1 
\end{lstlisting}
\caption{Numerical variation with Infinity Without Optimization (\texttt{-O0})}
\label{fig:case_study_2}
\end{figure}

Figure~\ref{fig:case_study_2} shows the code \& input-output 
generated by \varity, comprising three parts like the previous 
case study. The first part shows the \textbf{compute} method 
that performs a division operation involving the \texttt{ceil} function. 
The code was tested using both NVIDIA and AMD compilers to compare the results.

In the second part, we see a randomly generated discrepancy-inducing 
input, \texttt{+1.2374E-306}. With this input, the \texttt{nvcc} 
compiler produced an output of \texttt{Inf}, while the \texttt{hipcc} 
compiler produced \texttt{1.34887e-306}. To determine the cause of 
this variation, we analyzed the intermediate results and 
the assembly code generated by the compilers. 

The first statement of the \texttt{compute} method is an 
assignment statement and is the same for both devices. 
However, the discrepancy emerged when the \texttt{ceil} 
function was applied to \texttt{+1.5955E-125} in the following statement. 
The \texttt{nvcc} compiler computes \texttt{ceil(+1.5955E-125)} 
as 0, while the \texttt{hipcc} compiler results in 1. 
This difference in the \texttt{ceil} function's result 
led to a \textbf{division by zero} on the NVIDIA device, 
causing the \texttt{comp} value to become \texttt{Inf}. 
In contrast, the AMD device performed the division without issues, 
yielding \texttt{1.34887e-306}. The third part of 
Figure~\ref{fig:case_study_2} shows the difference in 
outputs for the same function call, \texttt{ceil(+1.5955E-125)}, across \texttt{nvcc} and \texttt{hipcc} compilers.

Further analysis revealed significant differences in 
how each platform implemented the \texttt{ceil} function.
NVIDIA's SASS and PTX assembly languages resulted in zero, 
while AMD's GPU ISA produced the expected result of 1. 

\subsubsection{Case Study 3: Inconsistencies in Infinity vs. NaN}

In this case study shown in Figure~\ref{fig:case_study_3}, 
the \texttt{compute} method performs several floating-point operations, 
such as calculating absolute values and using conditional statements. 
None of these operations cause any discrepancies, making this analysis 
particularly interesting. Initially, when no optimization 
flags are applied (\texttt{-O0}), both NVCC and HIPCC compilers 
produce consistent results, yielding \texttt{-inf} for the 
provided input values. This consistency suggests that, without optimizations, 
both compilers handle floating-point operations similarly.

\begin{figure}
\centering
\begin{lstlisting}[language=C]
__global__ /* __global__ is used for device run */
void compute(double comp, int var_1, double var_2, double var_3, double var_4, double var_5, double var_6, double var_7, double var_8) {
  double tmp_1 = (-1.8007E-323 - cosh(var_2 / -1.7569E192 + (-1.9894E-307 / +1.7323E-313 + var_3)));
  comp += tmp_1 + fabs(+1.5726E-307 - var_4);
  for (int i = 0; i < var_1; ++i) {
    comp += (+1.9903E306 / var_5);
  }
  if (comp >= (-1.4205E305 - (-1.4055E-312 * (var_6 + -1.7892E214 / var_7)))) {
    comp += +1.3803E305 * var_8;
  }
  printf("%.17g\n", comp);

}
\end{lstlisting}
\begin{lstlisting}[language=C, numbers=none]
Input: -1.5548E-320 5 +1.9121E306 +0.0 -1.1577E124 -1.8994E-311 +1.3675E306 +1.1296E-318 +1.2915E306

Output:
nvcc -O0:	-inf
hipcc -O0:	-inf
\end{lstlisting}
\begin{lstlisting}[language=C, numbers=none]
Input: -1.5548E-320 5 +1.9121E306 +0.0 -1.1577E124 -1.8994E-311 +1.3675E306 +1.1296E-318 +1.2915E306

Output:
nvcc -O1:	-inf
hipcc -O1:	-nan
\end{lstlisting}
\caption{Numerical variation with Infinity \& NaN Without Optimization (\texttt{-O0})}
\label{fig:case_study_3}
\end{figure}

However, discrepancies arise when any optimization 
flags are introduced. Notably, the \texttt{nvcc} compiler 
continues to produce \texttt{-inf}, while the \texttt{hipcc} 
compiler results in \texttt{-nan}. This divergence is 
not attributed to typical mathematical function behavior 
or rounding errors but rather to how optimizations impact 
the computation flow. Through detailed analysis, 
we found that discrepancies emerge 
after the optimization modifies 
intermediary computations within the code. 

The statement \texttt{comp += tmp\_1 + fabs(+1.5726E-307 - var\_4)} 
within the \texttt{compute} method produces \texttt{-inf} 
across all optimization levels and on both NVIDIA and AMD GPUs. 
This \texttt{-inf} value propagates through the subsequent operations, 
influencing the control flow when it encounters the \texttt{if} statement. 

As \texttt{comp} is compared within the branch, it remains \texttt{-inf}. 
However, once the execution exits the loop, optimizations potentially 
alter the sequence of calculations, leading to \texttt{comp} 
becoming \texttt{-nan}. This suggests that the optimization 
may adjust how intermediary values are handled, 
resulting in a propagation of \texttt{-nan} 
instead of \texttt{-inf}. 
The change from \texttt{-inf} to \texttt{-nan} 
in \texttt{hipcc} under optimization likely results 
from the reordering or elimination of intermediate steps, 
underscoring the complexities involved in maintaining 
consistency across different GPU architectures and compiler optimizations.

\begin{mdframed}[backgroundcolor=light-gray] 
\textbf{Answer to Q3:}
We identified the root causes of several inconsistencies. 
The root causes vary in each case. A key reason is differences 
in the low-level implementation of mathematical functions such as 
\texttt{fmod} and \texttt{ceil}. As demonstrated in Case 
Study 1, the \texttt{fmod} function was identified as a 
source of inconsistency. 
Compiler optimizations, particularly aggressive ones like \texttt{-ffast-math}, 
also contribute to discrepancies by affecting computation 
stability. Additionally, inconsistencies can arise 
from how intermediary values are handled due to different 
optimization strategies across compilers. 
\textbf{These factors collectively lead to variations in outputs on different GPU architectures.}
\end{mdframed}

\section{Related Work}
The \varity framework~\cite{varity-github,varity}, introduced by Laguna, utilizes random program generation to quantify floating-point variations in HPC systems, 
focusing on host-to-host and host-to-device testing. 
Our work extends this framework by incorporating device-to-device testing, 
offering a more detailed evaluation of numerical discrepancies 
between NVIDIA and AMD GPUs. Previous studies have examined the 
trade-offs between performance and numerical accuracy in compiler 
optimizations. Bentley et al.~\cite{bentley2018multi,bentley2019multi} 
conducted a multi-level analysis of compiler-induced variability while 
Guo et al.~\cite{guo2020pliner} developed PLiner to isolate lines of 
floating-point code contributing to variability. Chowdhary and Nagarakatte~\cite{chowdhary2021parallel} introduced parallel shadow 
execution to accelerate the debugging of numerical errors, 
highlighting the need for efficient tools in this domain.

Recent research has focused on detecting floating-point exceptions and 
inconsistencies across GPU architectures. Li et al.~\cite{li2023design} 
designed GPU-FPX, a tool for floating-point exception detection, 
and Innocente and Zimmermann~\cite{innocente2021accuracy} 
analyzed the accuracy of mathematical 
functions across different precision levels. 
Augonnet et al.~\cite{augonnet1} and Cao~\cite{cao1} 
explored task scheduling on heterogeneous architectures 
and the implications of using HIP technology on GPU-like accelerators, 
respectively, contributing to our understanding of how different architectures 
handle floating-point computations. 
In the context of heterogeneous systems, 
Miao et al.~\cite{miao2023expression} 
emphasized the importance of understanding compiler-induced inconsistencies. 
Additionally, Lopez-Novoa et al.~\cite{lopez2015} surveyed performance 
modeling techniques for accelerator-based computing, 
setting the stage for further exploration of discrepancies 
in numerical results between different GPU platforms. 
Finally, the adoption of alternative numerical formats 
like Posits has been explored as a means to enhance accuracy in HPC.
Poulos et al.~\cite{poulos2021} discussed the potential of 
Posits in improving precision in scientific applications, 
particularly in the context of exascale and edge computing. 
These studies underscore the ongoing challenges and 
opportunities in achieving consistent floating-point 
computations across different architectures and optimization levels.

\section{Limitations}
Our study has limitations. First, since programs are generated
randomly, we have no guarantee that these programs 
reflect the code in real scientific simulations. However,
the root cause of several cases arises from simple arithmetic and math
operations that are likely to be part of several larger codes.
Second, the \texttt{-ffast-math} option may not be directly comparable
between \texttt{nvcc} and \texttt{hipcc} as they most likely implement
optimizations in different ways; thus, programmers must be careful
when enabling them on code that runs on both GPUs.
Third, while most modern compilers use similar intermediate 
representations (e.g., LLVM IR), in general, optimization levels
are not necessarily comparable across compilers, i.e., \texttt{-O2}
in \texttt{nvcc} may not perform the same optimizations as 
\texttt{-O2} in \texttt{hipcc}. Nevertheless, it is crucial to
understand when numerical inconsistencies occur because programmers
often use the same optimization level when moving and running
code between different platforms.

\section{Conclusion}
Testing for numerical consistency between GPU platforms is crucial 
to ensure the correctness and robustness of numerical simulations.
As new HPC systems include different classes of GPUs, acceptance
testing is required to ensure system integrity and to fix issues
before production use---therefore testing for numerical consistency
between systems is essential.
This paper presents a study of compiler-induced numerical inconsistencies
between NVIDIA and AMD GPUs, which uses differential testing and 
random program generation.
We extended \varity to enable device-to-device testing, 
uncovering significant numerical discrepancies between 
NVIDIA and AMD GPUs. 
These discrepancies stem from various factors, 
including differences in the implementation of math functions
and discrepancies introduced by compiler optimizations, 
such as \texttt{-ffast-math}. Furthermore, based on our analysis
of assembly code, we conclude that
intermediary value handling differences across compilers seem to
contribute to these inconsistencies. 

In future work, we will investigate the root causes of 
mismatches between \HIPIFY-converted code and \varity-generated 
code. Additionally, we aim to develop automated debugging 
tools to efficiently identify and resolve these 
inconsistencies, minimizing manual analysis, and 
enhancing numerical consistency across heterogeneous 
computing platforms. 

{
\def\bibfont{\footnotesize} 
\bibliographystyle{ieeetr}
\bibliography{biblio}
}

\end{document}